%% file: main.tex
\documentclass[11pt,a4paper]{article}

\usepackage[T1]{fontenc}
\usepackage[utf8]{inputenc}
\usepackage[margin=1in]{geometry}
\usepackage[protrusion=true,expansion=false]{microtype}
\usepackage{amsmath,amssymb,amsthm,mathtools}
\usepackage{booktabs,tabularx,array}
\usepackage{graphicx}
\usepackage{enumitem}
\usepackage{xcolor}
\usepackage{url}
\usepackage[hidelinks]{hyperref}
\hypersetup{
  pdftitle={Hasse Refinement and Message-Passing Depth for Sheaf-Inspired Learning},
  pdfauthor={Neria Baris}
}

\newtheorem{theorem}{Theorem}[section]
\newtheorem{proposition}[theorem]{Proposition}
\newtheorem{lemma}[theorem]{Lemma}
\newtheorem{corollary}[theorem]{Corollary}
\theoremstyle{definition}
\newtheorem{definition}[theorem]{Definition}

\theoremstyle{remark}

\newcommand{\id}{\operatorname{id}}
\newcommand{\Face}{\operatorname{Face}}
\newcommand{\rank}{\operatorname{rk}}
\newcommand{\B}{\mathcal B}
\newcommand{\C}{\mathcal C}
\newcommand{\Hist}{\operatorname{Hist}}
\newcommand{\Hash}{\operatorname{HASH}}
\newcommand{\MSet}[1]{\left\{\!\left\{#1\right\}\!\right\}}
\newcommand{\R}{\mathbb R}
\newcommand{\GL}{\operatorname{GL}}

\setlist{leftmargin=*}
\title{Hasse Refinement and Message-Passing Depth\\for Sheaf-Inspired Learning}
\author{Neria Baris\\Mathematical Institute, University of Oxford}
\date{}

\begin{document}
\maketitle

\begin{abstract}
Cellular sheaves attach local data spaces to the cells of a combinatorial space, along with comparison maps across face relations. Motivated by sheaf-based learning on simplicial complexes, we study the distinguishing power and communication depth of local message passing. On finite graded posets, we compare color refinement on the rank-labeled Hasse graph with four-channel refinement through immediate faces, cofaces, and shared face or coface witnesses. If $h_t$ and $g_t$ are the respective colorings after $t$ rounds, then $h_{2t}\preceq g_t\preceq h_t$. The two rules induce the same stable partition on finite families. The factor-two gap in the first distinguishing round is attained by pairs of connected three-arm trees, viewed as one-dimensional simplicial complexes. Under bounded Hasse degree and exposed incidence labels, every $T$-layer four-channel message-passing network admits an exact functional simulation in $2T$ cover-only rounds, with explicit cache dimensions. A fixed cellular-sheaf Hodge operator is a linear specialization: its rank-$k$ Laplacian takes two rounds, and a degree-$p$ polynomial filter takes $2p$ rounds, a worst-case sharp bound under cell-local initialization. Experiments with factorable nonlinear messages recover the direct outputs to floating-point accuracy and measure the storage and execution costs of routing through witnesses.
\end{abstract}

\section{Introduction}

In recent years, the prominence of sheaf theory has been rising in the fields of data analysis and deep learning. This accelerated development began with the work of Hansen and Gebhart, who first introduced the architecture of sheaf neural networks as a generalization of Graph Convolutional Networks (GCNs), using fixed, predefined sheaves \cite{hansen2020sheaf}.

A cellular sheaf attaches local data spaces to the cells of a combinatorial space, along with comparison maps between local data across face relations. In a neural network, these maps specify how information carried by different cells is compared and combined. They allow the model to use context-dependent linear transformations when passing messages, accommodating signed or asymmetric relations and varying feature spaces \cite{hansen2020sheaf}. This extends the geometric deep learning framework to domains with nonuniform local comparisons \cite{bronstein2021geometric}.

Through this development, sheaves provide an elegant and formal way to represent the relationship between geometry and the data situated upon it, and between the local and the global. Sheaf-based models can also actively mitigate oversmoothing in deep learning on graphs \cite{bodnar2022neural}. Oversmoothing is the progressive loss of distinctions between node representations under repeated updates. By learning the comparison maps, these models can preserve patterns of variation that are compatible with their local relations, including in heterophilic settings \cite{bodnar2022neural}.

These advantages raise a question about the expressive power of local learning on sheaves: which distinctions can be made by its message-passing rules, and what depth is needed to make them? Weisfeiler--Leman color refinement provides a standard tool for comparing local message-passing architectures \cite{xu2019powerful}. On a simplicial complex, the relevant communication structure includes its higher-dimensional cells and their face relations. We therefore study refinement on the full face poset, where restriction maps and other incidence data can be exposed as labels.

Message passing on a simplicial complex often combines immediate faces and cofaces with same-dimensional cells that share a face or coface \cite{bodnar2021topological}. The latter interactions cross two cover relations, through the shared cell. Graded-poset refinement uses these four channels \cite{choi2026categorical}, whereas ordinary refinement on the Hasse graph uses one cover per round. We ask how this choice changes distinguishing power and the depth needed to reproduce an update. Both processes use the same lifted representation: each Hasse vertex represents a cell, so higher-order cells remain available to ordinary graph refinement.

Our main result is the round sandwich $h_{2t}\preceq g_t\preceq h_t$, for rank-refining initial labels. It gives equality of the stable partitions and bounds the gain in distinguishing depth by two. We prove that the bound is attained at object level for every $t$: a pair of connected three-arm trees is first distinguished after $t$ graded rounds and $2t$ Hasse rounds. Their maximum degree is three.

We then address the numerical question of reproducing a layer's output. With bounded Hasse degree and exposed incidence labels, we construct a target-aware Hasse simulator that reproduces a four-channel layer in two synchronous rounds. The labels, including any learned restriction maps, are available and held fixed during each simulated layer. Bounded-multiset encodings preserve arbitrary peer--witness messages, at an explicit cost in intermediate dimension; factorable messages use smaller sufficient statistics. An important linear example is the fixed cellular-sheaf Hodge Laplacian, whose action compares signals through adjacent ranks \cite{hansen2019spectral}. The incidence factorization computes it in two rounds, and Horner evaluation applies a degree-$p$ filter in $2p$ rounds. A path gives the matching worst-case lower bound. The numerical study compares explicit witnessed-peer records with their cover-only factorization using identical inputs and weights.

\section{Related Work}

\paragraph{Higher-order color refinement and simplicial message passing.}
Bodnar et al. introduced simplicial WL and message-passing simplicial networks with boundary, coboundary, lower, and upper neighborhoods \cite{bodnar2021topological}. Their analysis already contains factor-two round accounting when lower-adjacency information is routed through boundary incidences. Choi, Kim, and Yun formulate graded-poset WL, $F$-CatWL, and the corresponding four-message $F$-CatMPN \cite{choi2026categorical}. Their results compare these richer rules under channel reductions and injective operations. Here we instead compare the four-channel rule directly with ordinary one-hop refinement on the unaugmented undirected Hasse graph, which has no same-rank shortcut edges and no witness colors in its update.

\paragraph{Graph reductions of topological message passing.}
Several works encode higher-order message passing in graph representations. Hajij et al. use augmented Hasse representations \cite{hajij2023topological}, TopoTune uses neighborhood-specific augmented Hasse graphs \cite{papillon2025topotune}, and Taha et al. formulate simplicial and cellular message passing as relational structures \cite{taha2025demystifying}. Jogl et al. show that graph networks running on transformed representations can match their higher-order counterparts under the Weisfeiler--Leman criterion \cite{jogl2023expressivity}. Later, Chen et al. use the same four structural neighborhoods in CCWL and prove a bridge-neighborhood sufficiency result under coverage assumptions \cite{chen2026ccwl}. These constructions add relation-specific or same-rank connectivity. Here we retain only cover relations, so a mediated same-rank interaction must cross its witness and therefore takes two synchronous rounds.

\paragraph{Sheaf networks and Hodge operators.}
We use the covariant cellular convention and Hodge Laplacian of Hansen and Ghrist \cite{hansen2019spectral}. Learned sheaf diffusion and its relation to heterophily and oversmoothing were developed by Bodnar et al. \cite{bodnar2022neural}, and polynomial sheaf filters have also been studied directly \cite{borgi2025polynomial}. Bodnar et al. relate simplicial message passing to spectral convolutions of the Hodge Laplacian \cite{bodnar2021topological}. For arbitrary fixed finite-dimensional sheaf data, we ask how many communication rounds along cover relations are needed to apply the operator exactly, and what assumptions are needed for equivariance and sharpness.

Jogl et al.\ \cite[Theorem~F.1]{jogl2023expressivity} give a factor-two simulation on a transformed graph with one vertex per cell and upper-adjacency shortcut relations. Our bounds use only cover edges, include a matching connected-tree example, and give a functional numerical simulation with explicit cache dimensions. Table~\ref{tab:prior-comparison} compares the communication graphs and guarantees.

\begin{table}[ht]
\centering
\footnotesize
\begin{tabularx}{\textwidth}{@{}p{0.22\textwidth}XX@{}}
\toprule
Result & Communication and initialization & Guarantee relevant here \\
\midrule
Bodnar et al.\ \cite{bodnar2021topological}, Theorem~6
& Boundary and witnessed upper channels on simplicial complexes; the SWL definition starts from a constant color.
& Removing coboundary and lower channels preserves stable distinguishing power. \\
Jogl et al.\ \cite{jogl2023expressivity}, Theorem~F.1
& A transformed graph with upper-adjacency shortcut relations, one vertex per cell, and dimension features.
& A factor-two simulation of simplicial and cellular color refinement on that representation. \\
Choi et al.\ \cite{choi2026categorical}, Definitions~2 and~6
& Four witnessed channels on graded posets; Definition~2 starts from a constant color.
& Graded refinement and its corresponding numerical message-passing framework. \\
This paper
& Only cover edges; rank-refining initial labels, or equivalent directed incidence types for the left bound.
& A round sandwich, equal stable partitions, sharpness on connected trees, and functional simulation with explicit cache dimensions. \\
\bottomrule
\end{tabularx}
\caption{Communication graphs, initialization, and guarantees in the closest refinement and simulation results.}
\label{tab:prior-comparison}
\end{table}

\paragraph{Injective multiset aggregation.}
The constructive simulator requires finite-dimensional encodings that distinguish bounded multisets. The correspondence between injective aggregation and 1-WL is established \cite{xu2019powerful}. We use polarized power sums to obtain an explicit encoding for bounded cardinality, relying on classical separation results for multisymmetric polynomials \cite{rydh2007minimal} and related continuous representation results \cite{tabaghi2024universal}.

\section{Background and Mathematical Framework}

\subsection{Local communication on the face poset}

\begin{definition}[Simplicial complex and face poset]
A finite abstract simplicial complex $K$ is a finite family of nonempty finite sets closed under taking nonempty subsets. Its elements are called simplices or cells, and the dimension of a cell $\sigma$ is $|\sigma|-1$. Its face poset, denoted by $\Face(K)$, is the set of cells with the order relation
\[
             \sigma\leq\tau\quad\Longleftrightarrow\quad\sigma\subseteq\tau.
\]
The empty face is omitted from the poset used here.
\end{definition}

The face poset retains information absent from the graph skeleton alone. For example, the boundary of a triangle and a filled triangle have the same one-dimensional skeleton, with three vertices and three edges. The filled triangle also has a two-dimensional cell. Refinement that observes only the skeleton cannot detect this difference, whereas rank labels on the face poset detect it already in the initial histogram. Below we compare two refinement rules on the same face poset.

\begin{definition}[Graded poset and cover relation]
In a poset $P$, the relation $p\lessdot q$ holds when $p<q$ and there is no element $w$ with $p<w<q$. A poset is graded if it has a rank function $\rank:P\to\mathbb N$ satisfying
\[
                 \rank(q)=\rank(p)+1\qquad\text{when }p\lessdot q.
\]
In the face poset, the natural rank is dimension, and a cover relation is obtained by adding one vertex to a simplex.
\end{definition}

Each cell $\sigma$ carries a label containing its rank and original features,
\[
             \lambda(\sigma)=(\dim\sigma,x_\sigma).
\]
Write $P=\Face(K)$ and define the sets of immediate faces and cofaces by
\[
 \B(p)=\{w:w\lessdot p\},\qquad
 \C(p)=\{w:p\lessdot w\}.
\]
The Hasse graph is the undirected graph whose edges are cover relations. Rank in the labels determines which endpoint of an edge is the face and which is the coface. The labeled Hasse graph therefore determines the whole face poset by taking the transitive closure of the cover relations.

Write $\MSet{\cdot}$ for multisets and let $\Hash_H,\Hash_G$ be injective on the signatures under consideration.

Hasse--WL uses a single multiset of neighbors in the Hasse graph,
\[
 h_{r+1}(p)=\Hash_H\left(h_r(p),
    \MSet{h_r(q):q\in\B(p)\cup\C(p)}\right).
\]
The graded rule additionally observes same-rank peers through lower and upper witnesses. For a coloring $c$, define
\begin{align*}
 L_c(p)&=\MSet{(c(q),c(w)):
        w\lessdot p,\ w\lessdot q,\ q\ne p},\\
 U_c(p)&=\MSet{(c(q),c(w)):
        p\lessdot w,\ q\lessdot w,\ q\ne p}.
\end{align*}
The graded update is
\[
 g_{r+1}(p)=\Hash_G\left(
 g_r(p),
 \MSet{g_r(q):q\in\B(p)},
 \MSet{g_r(q):q\in\C(p)},
 L_{g_r}(p),U_{g_r}(p)
 \right).
\]
Each peer--witness path contributes one occurrence. To include incidence labels, assign to every cover $w\lessdot p$ a value $\ell_{wp}$ available in both communication directions, and write $\ell_{pw}=\ell_{wp}$. Define
\[
h_{r+1}^{\ell}(p)=\Hash_H\left(h_r^{\ell}(p),
\MSet{(h_r^{\ell}(q),\ell_{pq}):q\in\B(p)\cup\C(p)}\right).
\]
The labeled graded update replaces each direct record by $(g_r^\ell(q),\ell_{pq})$ and uses
\begin{align*}
L_{g_r^\ell}^{\ell}(p)&=\MSet{(g_r^\ell(q),g_r^\ell(w),\ell_{wp},\ell_{wq}):
w\lessdot p,q,\ q\ne p},\\
U_{g_r^\ell}^{\ell}(p)&=\MSet{(g_r^\ell(q),g_r^\ell(w),\ell_{pw},\ell_{qw}):
p,q\lessdot w,\ q\ne p}.
\end{align*}
For directed incidence types, append to each direct record its direction relative to the receiver; the same direction information is retained along each mediated path. The two processes use the same incidence values.

For two colorings on a common domain, write $a\preceq b$ if equality under $a$ implies equality under $b$, so that $a$ is the finer coloring. For object-level comparison we use the global histogram
\[
       \Hist_a(K)=\MSet{a(p):p\in\Face(K)}
\]
with one shared hash dictionary for all compared objects. For a graded poset $P$, use the same definition with $p\in P$.

\begin{definition}[Local synchronous communication model]
\label{def:local-model}
At round $r$, a cell $p$ has a finite-dimensional state $s_p^r$, including its rank or equivalent incidence-direction information. For adjacent cells, the message from $q$ to $p$ is
\[
 m_{q\to p}^r=M_r(s_p^r,s_q^r,\ell_{pq}),\qquad
 s_p^{r+1}=U_r\left(s_p^r,\MSet{m_{q\to p}^r:q\sim_Hp}\right).
\]
The update is permutation invariant in its multiset argument; messages and updates may be rank dependent, but do not depend on cell identifiers. Both endpoint states are read before any state is updated. Each message crosses one cover edge. Initialization depends only on the input at that cell and its exposed labels; it does not contain other cells' signals, a global state, or larger-radius features. Incidence labels $\ell_{pq}$ read in a block are available to both endpoint computations. They include every map, sign, and metric used by that block and are held fixed during it. For varying stalk dimensions, dimension and metric data are exposed locally and padding is fixed consistently. Local computation and the supplied labels have no communication-round charge; discovering or learning those labels is not simulated by this assumption.
\end{definition}

For a four-channel layer, the same pre-update convention applies. Write $z_p$ for the input state and $\ell_{wp},\ell_{wq}$ for labels on lower incidences. A lower-channel message has the form
\[
 f_{\downarrow}(z_p,z_q,z_w,\ell_{wp},\ell_{wq}),
 \qquad w\lessdot p,\ w\lessdot q,\ q\ne p,
\]
with an analogous upper message and direct messages depending on the target, neighbor, and incidence label. Each channel is aggregated separately by a permutation-invariant function before a local update. This defines the target-aware, edge-labeled function class used below. The lower bounds use the same endpoint-aware model: target dependence does not enlarge a one-round light cone.

\subsection{Reconstructing mediated channels in two rounds}

A four-channel message-passing layer uses direct boundary and coboundary messages together with lower- and upper-mediated messages. A mediated message through $w$ may depend jointly on the target state, peer state, witness state, and both incidence labels. A Hasse simulator therefore cannot replace a two-step interaction by a witness message that ignores the target.

To preserve this dependence, the first Hasse round stores an injective representation of the local neighborhood separated by rank. For a multiset of at most $m$ vectors in $\R^q$, define
\[
 J(q,m)=\binom{q+m}{q},\qquad
 \phi_{q,m}(x)=\bigl(1,(x^\alpha)_{1\leq|\alpha|\leq m}\bigr).
\]
Here $\alpha$ is a multi-index, $|\alpha|$ is the sum of its entries, and $x^\alpha$ is the corresponding monomial. A multiset $X$ is encoded by $E(X)=\sum_{x\in X}\phi_{q,m}(x)$. The first coordinate counts its elements, while the others are sums of monomials through degree $m$. This encoding is injective on multisets of cardinality at most $m$. A witness can therefore recover the states and edge labels in its neighborhood before the round.

In the second round, a target-aware message from the witness decodes the relevant neighborhood records, removes one occurrence of the target record, evaluates the original mediated message function on every remaining peer, and sends an injective encoding of the resulting message multiset. The receiver reconstructs the four channel multisets and applies the original aggregators and update. Duplicate states and labels are retained, since removing either of two identical target occurrences leaves the required multiplicity unchanged.

The linear case is simpler. If a mediated layer factors as
\[
 y_p=\sum_{w\sim_Hp}A_{pw}
       \left(\sum_{q\sim_Hw}B_{wq}x_q\right),
\]
then the first round computes $s_w=\sum_q B_{wq}x_q$ at the witness, and the second returns $A_{pw}s_w$ to the target. No multiset encoding has to be inverted.

\subsection{Cellular sheaves and local Hodge computation}

View $\Face(K)$ as a category with one morphism $\sigma\to\tau$ for each $\sigma\leq\tau$.

\begin{definition}[Cellular sheaf]
A real cellular sheaf on $K$ is a covariant functor
\[
             \mathcal F:\Face(K)\longrightarrow\mathbf{Vect}_{\R}.
\]
Explicitly, it assigns to each cell $\sigma$ a vector space $\mathcal F(\sigma)$, called its stalk, and to each face relation $\sigma\leq\tau$ a linear map
\[
             R_{\sigma\tau}:\mathcal F(\sigma)\longrightarrow\mathcal F(\tau).
\]
The maps satisfy $R_{\sigma\sigma}=\id_{\mathcal F(\sigma)}$ and
\[
             R_{\sigma\upsilon}=R_{\tau\upsilon}R_{\sigma\tau}
             \qquad(\sigma\leq\tau\leq\upsilon).
\]
\end{definition}

We use the covariant cellular convention, with maps from faces to cofaces \cite{hansen2019spectral}. The compatibility condition makes their composition independent of the chosen face chain. Simplex orientations supply incidence signs.

Assume all stalks are finite-dimensional and choose inner products on them. Choose an orientation on each simplex, and let $[\tau:\sigma]\in\{\pm1\}$ be the incidence number for a cover $\sigma^k\lessdot\tau^{k+1}$. The rank-$k$ cochain space is
\[
             C^k(K;\mathcal F)=\bigoplus_{\sigma\in K_k}\mathcal F(\sigma),
             \qquad K_k=\{\sigma\in K:\dim\sigma=k\}.
\]
The notation $R_{\sigma\tau}^*$ denotes the adjoint with respect to the chosen inner products.

The coboundary operator is
\[
 (\delta_kx)_\tau=\sum_{\sigma\lessdot\tau}
 [\tau:\sigma]R_{\sigma\tau}x_\sigma,
\]
and the rank-$k$ Hodge Laplacian is
\[
                  \Delta_k=
 \delta_{k-1}\delta_{k-1}^*+\delta_k^*\delta_k.
\]

A global section is a choice of $s_\sigma\in\mathcal F(\sigma)$ at every cell satisfying $R_{\sigma\tau}s_\sigma=s_\tau$ for every face relation. Thus a global section is a collection of local data consistent with all restriction maps. In rank zero, $\delta_0$ measures disagreement between the vertex data after mapping them into the common edge stalks, and $\ker\Delta_0=\ker\delta_0$ identifies with the space of global sections. More generally, finite-dimensional Hodge theory gives a harmonic representative of each sheaf cohomology class \cite{hansen2019spectral},
\[
 \ker\Delta_k\cong H^k(K;\mathcal F)
 =\ker\delta_k/\operatorname{im}\delta_{k-1}
\]
Applying the Hodge operator through local messages therefore implements a step in sheaf diffusion or spectral filtering. The two-round construction below gives its exact communication cost for fixed sheaf data.

We use zero operators outside the range of ranks, in particular $\delta_{-1}=0$. The composition condition and incidence signs give $\delta_{k+1}\delta_k=0$. The two-round Laplacian factorization below also holds for a family of cover maps that does not satisfy the composition condition. In that case it need not define a sheaf or a cochain complex.

Each of the two Laplacian summands factors through cells of an adjacent rank. Starting from a rank-$k$ signal $X_\sigma\in\mathcal F(\sigma)\otimes\R^c$, the first round stores
\begin{align*}
  L_\rho&=\sum_{\rho\lessdot\sigma}
       [\sigma:\rho]R_{\rho\sigma}^*X_\sigma,\\
  U_\tau&=\sum_{\sigma\lessdot\tau}
       [\tau:\sigma]R_{\sigma\tau}X_\sigma.
\end{align*}
The second round returns
\begin{align*}
  Y_\sigma^-&=\sum_{\rho\lessdot\sigma}
       [\sigma:\rho]R_{\rho\sigma}L_\rho,\\
  Y_\sigma^+&=\sum_{\sigma\lessdot\tau}
       [\tau:\sigma]R_{\sigma\tau}^*U_\tau,
\end{align*}
and outputs $Y_\sigma=Y_\sigma^-+Y_\sigma^+$. This computes $\Delta_kX$. Adding $X$ remains a separate local operation.

Messages use the incidence sign, restriction map, and stalk metrics to compute the adjoint. The maps and metrics are fixed during each two-round block, as in Definition~\ref{def:local-model}.

\section{Results and Analysis}

\subsection{Round-by-round distinguishing power}

We first compare the two refinement processes at cell level and at complex level.

\begin{theorem}[Round sandwich]
\label{thm:round-sandwich}
Let $P$ be a finite graded poset and let $h_0=g_0=\lambda$, where $\lambda$ is rank-refining. Then, for every $t\geq0$,
\[
                         h_{2t}\preceq g_t\preceq h_t.
\]
The right refinement holds for arbitrary common initial labels. For the left refinement it is enough that the direction of every Hasse incidence is recoverable, either from rank or from an explicit edge type recording its direction. The statement remains valid with incidence labels when the same labels are exposed to both processes.
\end{theorem}

The direct boundary and coboundary multisets in a graded update determine their union, which is exactly the one-hop Hasse neighborhood, giving $g_t\preceq h_t$. In the other direction, two Hasse rounds expose the labeled neighborhood of each immediate witness. Removing one occurrence of the target record from each decoded neighborhood reconstructs the same-rank peers reached through that witness. Aggregating over witnesses reconstructs both mediated channels with the correct multiplicities. Induction therefore gives $h_{2t}\preceq g_t$.

The cell-level comparison passes immediately to objects through global histograms.

\begin{corollary}[Object-level consequences]
\label{cor:object-sandwich}
At round $t$, every object distinction made by $h_t$ is made by $g_t$, and every object distinction made by $g_t$ is made by $h_{2t}$. On a finite family, the stable cell partitions and the object partitions induced by their global histograms coincide.
\end{corollary}

\subsection{Sharpness of the factor-two separation}

We now ask whether the factor two can be realized at object level without isolated cells or disconnected unions. For positive integers $a_1,a_2,a_3$, let $S(a_1,a_2,a_3)$ be the one-dimensional complex of the tree with one degree-three center and three internally disjoint center-to-leaf paths containing $a_1,a_2,a_3$ graph edges, respectively.

\begin{theorem}[Connected object-level sharpness]
\label{thm:connected-sharpness}
For every $t,M\geq1$, put
\[
 K_{t,M}=S(2t-1,2t+1,M),\qquad
 L_{t,M}=S(2t,2t,M).
\]
With rank initialization and joint global hashes,
\begin{align*}
 \Hist_{h_r}(K_{t,M})&=\Hist_{h_r}(L_{t,M}) &&(0\leq r<2t),\\
 \Hist_{h_{2t}}(K_{t,M})&\ne\Hist_{h_{2t}}(L_{t,M}),\\
 \Hist_{g_s}(K_{t,M})&=\Hist_{g_s}(L_{t,M}) &&(0\leq s<t),\\
 \Hist_{g_t}(K_{t,M})&\ne\Hist_{g_t}(L_{t,M}).
\end{align*}
The complexes are nonisomorphic connected trees. As one-dimensional complexes they are pure and flag. Both their maximum graph degree and the maximum degree of their Hasse graphs are three.
\end{theorem}

Subdividing each graph edge once turns the Hasse graph into a three-arm tree with alternating ranks, so an arm with $a$ graph edges has Hasse length $2a$. Before the endpoint waves meet, the color bands depend only on distance from the center or leaf, while the interior cell counts depend affinely on total arm length. The two varying arms have the same total length in $K_{t,M}$ and $L_{t,M}$. The first Hasse discrepancy occurs when the two endpoint profiles meet at round $2t$. The mediated graded rule passes the corresponding information between same-rank cells through a witness, reaching the meeting at round $t$.

\subsection{Exact functional simulation of four-channel layers}

\begin{theorem}[Exact two-round simulation of each layer]
\label{thm:catmpn-simulation}
Assume every cell satisfies $|\B(p)|\leq b$ and $|\C(p)|\leq c$. Assume also that every incidence value read at layer $t$ is available as an edge label of dimension at most $e_t$. Then every $T$-layer four-channel CatMPN in the target-aware, edge-labeled message-function class has a Hasse--MPN realization of depth $2T$. After every even Hasse round, the state coordinates simulating the original layer equal the CatMPN states numerically.

Let $d_t$ be the cell-state dimension after padding to a common dimension, and let $R$ be the maximum rank. Put $q_t=d_t+e_t$ and
\[
 L=b\max\{c-1,0\},\qquad U=c\max\{b-1,0\}.
\]
If the four message dimensions are $a_{B,t},a_{C,t},a_{\downarrow,t},a_{\uparrow,t}$, sufficient cache and second-round widths are
\begin{align*}
 W_t^{\mathrm{cache}}
 &= (R+1)+d_t+J(q_t,b)+J(q_t,c),\\
 P_t^{(2)}
 &=J(a_{B,t},b)+J(a_{C,t},c)
   +J(a_{\downarrow,t},L)+J(a_{\uparrow,t},U).
\end{align*}
\end{theorem}

The proof uses the two-round construction of Section~3.2. Round one stores separate encodings of the lower and upper neighborhoods. Round two reconstructs the direct and mediated message multisets, after which the receiver applies the original aggregators and update. In this construction, the cache for the next layer is built only after the state $z^{t+1}$ is available, since a synchronous message is computed from the pre-update state. Repeating the block gives $2T$ rounds. The construction gives a network-level upper bound; the one-round obstruction below supplies the layer-level worst case.

Exactness here is equality of the composed functions. On a fixed finite realized domain, finite interpolation realizes the required maps. On compact state domains with continuous layer components, the encoders have continuous inverses on their images. Continuous extensions of the decoder components then give uniform MLP approximation through the finite sequence of simulation blocks.

The general-purpose cache can be large: $J(64,10)=718{,}406{,}958{,}841$. Proposition~\ref{prop:encoding-width} gives the classical dimension bound $mq$ for continuous injective encodings of multisets of $m$ vectors in an open subset of $\mathbb R^q$. Proposition~\ref{prop:moment-conditioning} exhibits a non-Lipschitz moment inverse at repeated records. Linear and factorable layers instead admit compact sufficient-statistic caches.

The two-round upper bound is also necessary in the worst case. Under rank-only local initialization, an isolated edge and a distinguished edge of a two-edge path have the same radius-one Hasse neighborhood. A four-channel layer can nevertheless distinguish them by counting lower-mediated peers. Hence no anonymous one-round synchronous Hasse--MPN can simulate every full four-channel layer in the stated local model.

\subsection{Sheaf--Hodge layers and polynomial filters}

For sheaf-based learning, the face poset supplies the communication graph, while the restriction maps determine how local data are compared. The Hodge realization applies to a fixed sheaf within each operator block. The general simulator also covers nonlinear messages that use the same exposed incidence maps, with its cache dimensions measuring the cost of reproducing witnessed interactions through cover relations.

The general cache is unnecessary when the target layer is the fixed Hodge operator described in Section~3.3.

\begin{theorem}[Exact and sharp Hasse realization]
\label{thm:sheaf-laplacian}
Fix a rank $k$ and a finite oriented simplicial complex with finite-dimensional real inner-product stalks and fixed restriction maps. Under the exposed-incidence and cell-local initialization assumptions of Section~3.3, two synchronous Hasse rounds, in which every message crosses one edge, compute exactly $\Delta_kX$. If $X$ or $XA$ is retained locally, the same two rounds followed by a local update compute
\[
              X\longmapsto\Phi\bigl(XA+(\Delta_kX)B\bigr).
\]
Here $A,B$ act on feature channels and $\Phi$ is a local function. If restriction maps and stalk metrics remain fixed, a degree-$p$ polynomial filter with matrix coefficients on the feature channels is computed in $2p$ rounds by the Horner recurrence
\[
       H_p=XA_p,\qquad
       H_j=\Delta_kH_{j+1}+XA_j\quad(j=p-1,\ldots,0).
\]
For a scalar polynomial $q(z)=\sum_{j=0}^pa_jz^j$ with $p\geq1$ and $a_p\ne0$, the bound $2p$ is worst-case sharp already for $k=0$, in the model with cell-local initialization and communication only along cover relations.
\end{theorem}

At every stage, the maps land in the vector space in which the messages must be summed. Downward first-round messages $R_{\rho\sigma}^*X_\sigma$ lie in $\mathcal F(\rho)$, upward messages $R_{\sigma\tau}X_\sigma$ lie in $\mathcal F(\tau)$, and both return maps land in $\mathcal F(\sigma)$. The diagonal blocks are return-to-source length-two paths, and remain part of the Laplacian.

The sharpness argument uses the constant rank-one sheaf on a path. For vertices separated by $p$ graph edges, the Hasse distance between the endpoints is $2p$, while the leading term of $q(\Delta_0)$ produces a nonzero endpoint-to-endpoint coefficient. In fewer than $2p$ synchronous rounds along cover relations, the endpoint state is independent of the initial signal at the opposite endpoint. Comparing the zero input with a unit impulse therefore gives an all-input operator lower bound.

The linear computation of $\Delta_kX$ is equivariant under orthogonal changes of basis in each stalk. This does not automatically extend to a layer containing $\Phi$. Such a layer additionally requires $\Phi(QX)=Q\Phi(X)$ for every permitted change of basis $Q$. For example, coordinatewise ReLU is not equivariant under all orthogonal basis changes. More generally, an arbitrary basis change $G_\sigma\in\GL(\mathcal F(\sigma))$ is compatible with the same coordinate computation only when stalk metrics are transported as well and adjoints are taken in the transported metrics. With fixed Euclidean metrics and transposes that ignore the metric, general $\GL$ equivariance fails.

Changing the sheaf maps between the two rounds changes the operator. If round one uses $R^{(1)}$ and round two uses $R^{(2)}$, the output is
\[
 \delta_{k-1}^{(2)}(\delta_{k-1}^{(1)})^*X
 + (\delta_k^{(2)})^*\delta_k^{(1)}X,
\]
which is generally neither self-adjoint nor positive semidefinite. Likewise, changing the Laplacian between Horner blocks produces ordered products of different operators rather than a polynomial in one fixed Laplacian.

\input{experiment_section.tex}

\begin{figure}[ht]
\centering
\includegraphics[width=\textwidth]{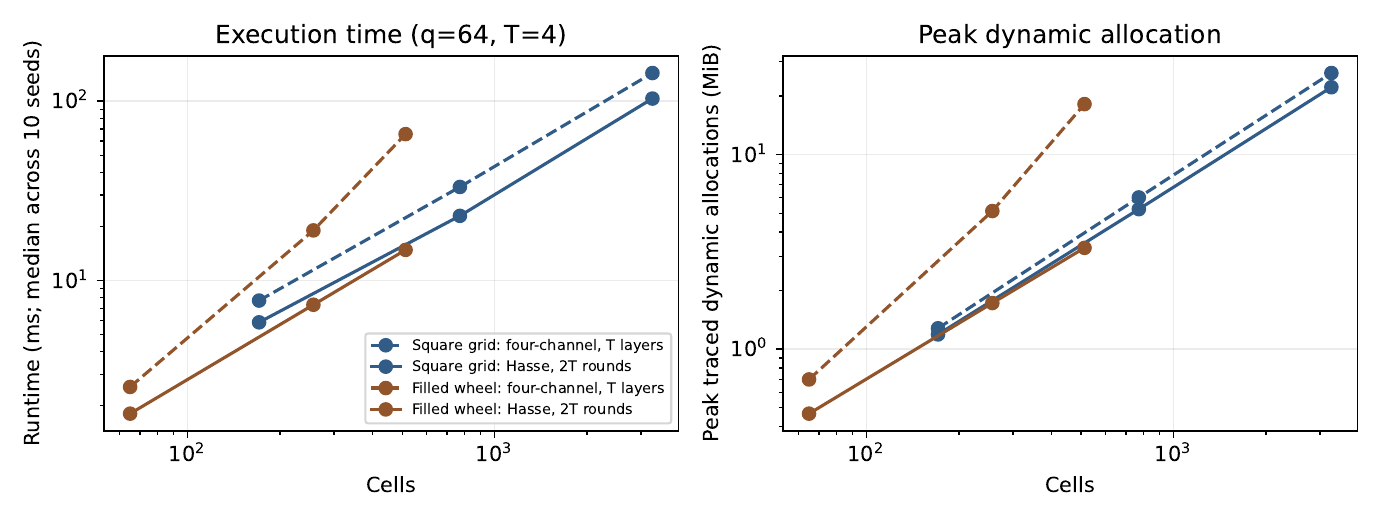}
\caption{Execution time and peak dynamic allocation at $q=64,T=4$. Dashed curves use explicit witnessed-peer records; solid curves use cover-only routing. Allocation during execution excludes prebuilt topology, inputs and weights.}
\label{fig:depth-resources}
\end{figure}

\section{Discussion}

The distinction between a witnessed update and cover-only communication is operational: one computes a two-incidence interaction directly, while the other routes it through its witness. The general simulator retains enough information to evaluate arbitrary peer--witness messages. Factorable messages can compress that information, as in the nonlinear study and the sheaf-Laplacian specialization.

The bounds assume rank or direction information, cell-local initialization, and exposed incidence values fixed within each simulation block. Shortcut edges, global communication, or larger-radius initial features change the communication model. Exact functional simulation and continuous approximation have separate parameter and conditioning costs: finite ReLU maps can approximate the required continuous components, and the lossless-cache dimension bound applies to injective neighborhood encodings. The numerical study uses fixed weights and compares explicit peer materialization with factored sums; its runtime and allocation measurements concern those implementations. Prediction accuracy, training behavior, and optimal cache dimension remain separate questions.

Useful extensions include identifying message families with smaller sufficient statistics, obtaining approximation bounds that track decoder conditioning, and allowing learned incidence maps to change while retaining an operator interpretation.

\paragraph{Reproducibility.}
Code, configurations and raw measurements accompany the paper. The supplement documents verification, the numerical protocol and the exploratory denoising study.

\bigskip
\appendix

\section{Technical Setup and Proofs}

\subsection{The face functor and object comparison}

Let $\mathbf{Simp}_{\mathrm{lab}}$ be the category of finite labeled simplicial complexes with label-preserving simplicial isomorphisms, and let $\mathbf{GPos}_{\mathrm{lab}}$ be the category of finite labeled graded posets with label-preserving, rank-preserving order isomorphisms. The face functor
\[
 \Face:\mathbf{Simp}_{\mathrm{lab}}\longrightarrow\mathbf{GPos}_{\mathrm{lab}},\qquad
 \Face(K,x)=(\Face(K),\lambda),
\]
maps a vertex bijection to the induced map on faces. It preserves identities and composition, is faithful on these categories, and reflects object isomorphism because rank-zero elements recover the vertices and each face is determined by the rank-zero elements below it.

All compared objects are colored on their disjoint union, or equivalently with one canonical global hash dictionary. Since rank remains part of every color, the all-cell histogram is equivalent to the tuple of rank-wise histograms.

\subsection{Proof of the round sandwich}

\begin{proof}[Proof of Theorem~\ref{thm:round-sandwich}]
We first prove $g_t\preceq h_t$. The case $t=0$ is equality. Assume the claim at $t$, and suppose $g_{t+1}(p)=g_{t+1}(q)$. Injectivity of $\Hash_G$ gives equal central colors and equal direct lower and upper multisets. By the induction hypothesis, replacing each $g_t$ color by the $h_t$ color it determines preserves equality of these multisets. Their union is the Hasse neighbor multiset, so the two $h_{t+1}$ signatures agree. This direction uses no rank information.

We now prove $h_{2t}\preceq g_t$. Again the base case is equality. Assume $h_{2t}\preceq g_t$ and suppose $h_{2t+2}(p)=h_{2t+2}(q)$. Decoding the two outer Hasse hashes gives identical old central $h_{2t}$ colors and equal multisets, separated by neighbor rank, of
\[
 \left(h_{2t}(w),\MSet{h_{2t}(z):z\sim_Hw}\right)
\]
for witnesses adjacent to $p$ and $q$. For a lower witness $w\lessdot p$, retain only neighbors of rank $\rank p$ in its decoded neighborhood, and delete one occurrence of $h_{2t}(p)$. The remainder is exactly the multiset of same-rank peers reached through $w$. Deleting a single occurrence is well-defined even if another peer has the same color, since the required multiplicity is unchanged. Pair each remaining peer color with the witness color, and map both to $g_t$ colors using the induction hypothesis. Aggregating over lower witnesses reconstructs $L_{g_t}(p)$, and the same argument reconstructs $U_{g_t}(p)$ from upper witnesses. The cell's decoded neighborhood also determines the direct boundary and coboundary multisets. Thus every input to the next graded update agrees at $p$ and $q$, giving $g_{t+1}(p)=g_{t+1}(q)$.

If directions are exposed as edge types rather than ranks, the same separation argument applies. For incidence labels, a decoded witness star contains pairs of a color and an incidence value. Delete one occurrence of $(h_{2t}^{\ell}(p),\ell_{wp})$ for a lower witness, using the outer incidence record to identify the value. The upper case is symmetric. Each remaining record retains its peer incidence value, so both values in the labeled mediated channel are reconstructed.
\end{proof}

\begin{proof}[Proof of Corollary~\ref{cor:object-sandwich}]
If $a\preceq b$, each realized $a$ color determines a unique $b$ color. Applying this map to the $a$ histogram gives the corresponding $b$ histogram. The two object-level claims follow from the two refinement relations in Theorem~\ref{thm:round-sandwich}. On a finite family, choose $t$ beyond both stabilization times.
\end{proof}

\subsection{Proof of the constructive simulator}

For a multiset $X$ of at most $m$ vectors in $\R^q$, the cardinality coordinate and polarized power sums through degree $m$ give an injective encoding $\sum_{x\in X}\phi_{q,m}(x)$. In characteristic zero, for a fixed cardinality $n\leq m$, power sums through degree $n$ generate the algebra of multisymmetric polynomials on $n$ vectors and separate permutation orbits. The cardinality coordinate separates different sizes \cite{rydh2007minimal}.

\begin{proof}[Proof of Theorem~\ref{thm:catmpn-simulation}]
Fix a layer $t$. In the first Hasse round, lower and upper neighbors write to separate blocks, and every cell stores
\begin{align*}
 E_p^-&=\sum_{w\in\B(p)}\phi_{q_t,b}(z_w^t,e_{wp}),\\
 E_p^+&=\sum_{w\in\C(p)}\phi_{q_t,c}(z_w^t,e_{pw}),
\end{align*}
together with its rank and old state $z_p^t$. These sums injectively encode the rank-separated Hasse neighborhoods. Together with the $d_t$ state coordinates and $R+1$ coordinates for rank encoding, they give the stated cache width. A first-round message width of $J(q_t,b)+J(q_t,c)$ and block output width of $(R+1)+d_{t+1}$ suffice.

Consider a second-round message from a lower witness $w\lessdot p$. The sender decodes its stored upper neighborhood, while the target-aware message supplies $(z_p^t,e_{wp})$. Removing one occurrence of that pair leaves exactly the peer records needed by the lower-mediated channel. The message evaluates the original mediated message function on each remaining peer and emits the moment encoding of the resulting multiset. It also emits the moment encoding of the original direct boundary message. An upper witness performs the symmetric construction for the direct coboundary and upper-mediated channels.

After summing the incoming messages, the receiver's four blocks encode exactly the four message multisets of the original layer. Their cardinalities are bounded by $b$, $c$, $L$, and $U$ as defined above, giving the displayed width bound. These bounds also cover empty mediated channels. The receiver decodes the multisets, applies the original permutation-invariant aggregators, and then the original update. Its output coordinates simulating the original layer are therefore $z_p^{t+1}$.

Since a synchronous message is computed from the pre-update state, the cache for layer $t+1$ cannot be built in the same round that produces $z^{t+1}$. Repeating the two-round block for $t=0,\ldots,T-1$ gives exact equality after every even Hasse round and total depth $2T$.
\end{proof}

\begin{proposition}[Linear witnessed specialization]
\label{prop:linear-witness}
If
\[
 y_p=\sum_{w\sim_Hp}A_{pw}\left(\sum_{q\sim_Hw}B_{wq}x_q\right),
\]
then two Hasse rounds compute the layer without an injective cache. If a diagonal path is excluded, its contribution $A_{pw}B_{wp}x_p$ can be subtracted locally.
\end{proposition}

\begin{proof}
Round one stores $s_w=\sum_qB_{wq}x_q$ at every witness. Round two sends $A_{pw}s_w$ to $p$ and sums over witnesses, giving the displayed double sum. If the path $p\to w\to p$ is forbidden, its contribution is already known at the target and can be subtracted. No multiset encoding has to be inverted.
\end{proof}

\begin{proposition}[Worst-case one-round obstruction]
\label{prop:one-round-obstruction}
Under rank-only local initialization, no anonymous one-round synchronous Hasse--MPN can simulate every full four-channel CatMPN at cell level when messages cross one Hasse edge, there is no global state or identifier, and messages are computed from pre-round states.
\end{proposition}

\begin{proof}
An isolated edge and a distinguished edge of a two-edge path have the same rank-labeled radius-one Hasse neighborhood. A four-channel layer returning the number of lower-mediated edge peers gives zero in the first case and one in the second. A one-round local Hasse process therefore cannot match both outputs.
\end{proof}

\subsection{Cache width and numerical conditioning}

\begin{proposition}[A lower bound for continuous lossless encodings]
\label{prop:encoding-width}
Let $q,m\geq1$ and let $V\subset\mathbb R^q$ be nonempty and open. If a continuous permutation-invariant encoder of all size-$m$ multisets in $V$ is injective and takes values in $\mathbb R^s$, then $s\geq mq$. The same bound holds for an encoder of multisets of cardinality at most $m$.
\end{proposition}

\begin{proof}
Choose $m$ pairwise disjoint open balls $B_1,\ldots,B_m$ in $V$. On their product the unordered encoding is injective even as a function of ordered tuples: an element in $B_i$ cannot be exchanged with one in another ball. If $s<mq$, compose this continuous injection from the open subset $\prod_i B_i\subset\mathbb R^{mq}$ with the inclusion $\mathbb R^s\hookrightarrow\mathbb R^{mq}$. Invariance of domain \cite[Theorem~2B.3]{hatcher2002algebraic} makes the image open, although it lies in a proper linear subspace, a contradiction. Restricting an at-most-$m$ encoder to size $m$ proves the last statement.
\end{proof}

\begin{proposition}[Moment decoding can be poorly conditioned]
\label{prop:moment-conditioning}
Equip two-element scalar multisets with the minimum Euclidean matching distance. The inverse of the moment encoding $E(X)=(|X|,\sum x,\sum x^2)$ is not Lipschitz on any neighborhood of $\{0,0\}$ containing $\{-a,a\}$ for arbitrarily small $a>0$.
\end{proposition}

\begin{proof}
For $X_a=\{-a,a\}$ and $X_0=\{0,0\}$, the matching distance is $\sqrt2a$, whereas $\|E(X_a)-E(X_0)\|_2=2a^2$. Their ratio $1/(\sqrt2a)$ is unbounded as $a\to0$.
\end{proof}

\subsection{Proof of the sheaf realization}

\begin{proof}[Proof of Theorem~\ref{thm:sheaf-laplacian}]
It suffices to check one feature channel. For $\rho^{k-1}\lessdot\sigma^k$,
\[
 R_{\rho\sigma}:\mathcal F(\rho)\to\mathcal F(\sigma),\qquad
 R_{\rho\sigma}^*: \mathcal F(\sigma)\to\mathcal F(\rho).
\]
Thus the downward first-round message lies in $\mathcal F(\rho)$, while $R_{\sigma\tau}X_\sigma$ lies in $\mathcal F(\tau)$ for $\sigma^k\lessdot\tau^{k+1}$. The return maps land in $\mathcal F(\sigma)$, so every sum combines messages in the same vector space. By definition,
\[
 L=\delta_{k-1}^*X,\qquad U=\delta_kX,\qquad
 Y^-=\delta_{k-1}L,\qquad Y^+=\delta_k^*U.
\]
Hence $Y=\Delta_kX$. Tensoring with the identity on feature channels gives the multichannel statement.

For cells $\sigma,\sigma'\in K_k$, the block from $\mathcal F(\sigma')$ to $\mathcal F(\sigma)$ is
\begin{align*}
 (\Delta_k)_{\sigma,\sigma'}={}&
 \sum_{\substack{\rho\in K_{k-1}\\
        \rho\lessdot\sigma,\ \rho\lessdot\sigma'}}
 [\sigma:\rho][\sigma':\rho]
 R_{\rho\sigma}R_{\rho\sigma'}^*\\
 &+\sum_{\substack{\tau\in K_{k+1}\\
        \sigma\lessdot\tau,\ \sigma'\lessdot\tau}}
 [\tau:\sigma][\tau:\sigma']
 R_{\sigma\tau}^*R_{\sigma'\tau}.
\end{align*}
The first composition has type $\mathcal F(\sigma')\to\mathcal F(\rho)\to\mathcal F(\sigma)$, and the second has type $\mathcal F(\sigma')\to\mathcal F(\tau)\to\mathcal F(\sigma)$. Setting $\sigma=\sigma'$ shows that the positive diagonal return paths remain included.

For the residual layer, retain $XA$ while the two communication rounds compute $\Delta_kX$, multiply the latter by the channel factor $B$, add, and only then apply $\Phi$. An intermediate activation changes an incidence factor and no longer computes the stated operator.

For the polynomial filter, descending induction gives
\[
 H_j=\sum_{r=0}^{p-j}\Delta_k^rXA_{j+r}.
\]
The identity holds for $j=p$. Applying $\Delta_k$ and adding $XA_j$ gives the next case, so $H_0$ is the required polynomial. There are $p$ applications of a fixed Laplacian, each taking two communication rounds.

For sharpness, take the path $v_0-e_1-v_1-\cdots-e_p-v_p$ with the constant rank-one sheaf and orient each edge from $v_{i-1}$ to $v_i$. Then $\Delta_0$ is the graph Laplacian. For $j<p$, the coefficient $(\Delta_0^j)_{v_p,v_0}$ is zero because $j$ factors cannot cross $p$ graph edges. For $j=p$, the unique monotone walk contributes $(-1)^p$, so
\[
             (q(\Delta_0))_{v_p,v_0}=a_p(-1)^p\ne0.
\]
The Hasse distance from $v_0$ to $v_p$ is $2p$. Under cell-local initialization and without a global or nonlocal channel, induction on synchronous rounds confines each state to its Hasse light cone. In fewer than $2p$ rounds, the state at $v_p$ is independent of the input at $v_0$. Comparing the zero input with a unit impulse at $v_0$ proves the all-input lower bound.

Let $Q_k=\bigoplus_{\sigma\in K_k}Q_\sigma$. Orthogonal gauges satisfy
\[
 \delta'_k=Q_{k+1}\delta_kQ_k^*,\qquad
 (\delta'_k)^*=Q_k\delta_k^*Q_{k+1}^*,\qquad
 \Delta'_k=Q_k\Delta_kQ_k^*.
\]
Substitution gives $L'_\rho=Q_\rho L_\rho$, $U'_\tau=Q_\tau U_\tau$, and $Y'_\sigma=Q_\sigma Y_\sigma$. For a general basis change $x'=Gx$, transport the metric as $M'=G^{-*}MG^{-1}$. If $R'=G_\tau RG_\sigma^{-1}$, its metric adjoint satisfies $(R')^{\dagger'}=G_\sigma R^\dagger G_\tau^{-1}$. The same round-by-round equivariance therefore holds with $G$ in place of $Q$. This statement concerns the linear computation. A nonlinear layer additionally requires an equivariant activation function.

For the raw Laplacian on $c$ channels, an active intermediate cell $\omega$ stores $\dim\mathcal F(\omega)c$ scalars. Horner evaluation stores the input and current accumulator at rank $k$, so local memory is linear in stalk and channel dimensions. Excluding restriction maps, a scalar polynomial has $p+1$ coefficients and a matrix-channel polynomial has $(p+1)c_{\rm in}c_{\rm out}$ coefficients, independent of the number of cells.
\end{proof}

\begin{proposition}[Maps changed inside a block]
\label{prop:mixed-operator}
If round one uses restriction maps $R^{(1)}$ and round two uses $R^{(2)}$, with fixed stalk metrics, the output is
\[
 \delta_{k-1}^{(2)}(\delta_{k-1}^{(1)})^*X
 + (\delta_k^{(2)})^*\delta_k^{(1)}X.
\]
This need not be the Hodge Laplacian of either restriction family.
\end{proposition}

\begin{proof}
Round one stores $(\delta_{k-1}^{(1)})^*X$ and $\delta_k^{(1)}X$. The return maps in round two are $\delta_{k-1}^{(2)}$ and $(\delta_k^{(2)})^*$, giving the displayed operator. Its adjoint interchanges the superscripts, so it is generally not self-adjoint. A Hodge Laplacian built from a single map family is self-adjoint and positive semidefinite.
\end{proof}

\subsection{Proof of connected sharpness}

\begin{lemma}[An auxiliary marked line-graph profile]
\label{lem:line-profile}
For a simple graph viewed as a one-dimensional complex, initialize line-graph refinement by $f_0(e)=\deg_{L(G)}e$ and let $f_{r+1}$ be ordinary color refinement on $L(G)$. Under rank initialization, $g_s$ on the edge cells refines $f_{s-1}$ for every $s\geq1$. In particular, $g_s(e)$ determines the degree-marked nonbacktracking unfolding of $L(G)$ to radius $s-1$. On an interior edge of a spider arm, stop each of its two rays at the first degree mark different from $2$; a stopped ray lying within that radius, including its endpoint mark and distance, is determined by the color.
\end{lemma}

\begin{proof}
At the first graded round, the lower-mediated multiplicity is $\deg u+\deg v-2=\deg_{L(G)}(uv)$, so $g_1$ determines $f_0$. At each later round, discard witness colors from the lower-mediated pairs. Because distinct edges in a simple graph share at most one vertex, the resulting multiset contains each line-graph neighbor once. Induction gives $g_s\preceq f_{s-1}$ on edge cells. Ordinary refinement retains the central color and the multiset of neighboring colors. Recursively decoding these signatures recovers the rooted nonbacktracking unfolding: at each extension remove one occurrence of the incoming parent's earlier color, retaining all other multiplicities. At radius zero the degree mark supplies the number of available continuations. Induction recovers the marked radius-$r$ unfolding from $f_r$. Before a mark different from $2$ is reached on an arm, every nonbacktracking continuation is unique; hence stopping there and recording its distance is determined by the unfolding.
\end{proof}

\begin{proof}[Proof of Theorem~\ref{thm:connected-sharpness}]
The center is the unique vertex of graph degree three, so the multiset of arm lengths is an isomorphism invariant. The multisets $\{2t-1,2t+1,M\}$ and $\{2t,2t,M\}$ cannot agree. Each complex is a tree without isolated vertices, hence connected, pure, one-dimensional, and flag, with maximum graph degree three. Subdividing each graph edge once identifies its Hasse graph with an alternating-rank three-arm tree, giving Hasse arm length $\ell=2a$ for $a$ graph edges.

The case $r=0$ follows from equal vertex and edge counts. Assume henceforth that $r\geq1$. On a rank-colored tree, $h_r(x)$ is determined by, and determines, the isomorphism type of the rooted radius-$r$ neighborhood, with initial rank labels at its boundary. A center or leaf at distance $d$ first changes the root color at round $d+1$. Before the two endpoint waves overlap, the noncenter positions $1,\ldots,\ell$ on an arm split into three regions,
\[
 \underbrace{1,\ldots,r-1}_{\text{center band}},\quad
 \underbrace{r,\ldots,\ell-r}_{\text{interior}},\quad
 \underbrace{\ell-r+1,\ldots,\ell}_{\text{leaf band}}.
\]
The color of each band is independent of total arm length, while rank-wise interior counts depend affinely on the even number $\ell$. For $r<2t$, the shortest varying Hasse arm has length $4t-2$, and $2(r-1)<4t-2$, so the bands do not overlap. Both objects have two varying arms with the same total Hasse length,
\[
 2(2t-1)+2(2t+1)=2(2t)+2(2t)=8t.
\]
From the center and every cell on the common $M$ arm, a varying-arm leaf is at distance at least $4t-2>r$. For a center-band cell at position $d\leq r-1$, its own varying-arm leaf is at distance at least $4t-2-(r-1)\geq2t>r$, and the other varying-arm leaf is farther away. Thus the center color agrees in both objects, as do all common-arm views and corresponding center-band views. A leaf-band view does not reach the center, so it depends only on distance from its leaf. Each band color therefore has the same multiplicity, and the rank-wise interior counts agree. This covers the center as well as all arm positions and gives equality of Hasse histograms through round $2t-1$.

At round $2t$, Hasse position $2t-1$ on the short graph arm of length $(2t-1)$ edges is an edge cell whose radius-$2t$ view sees the exceptional center and leaf degrees in opposite components. This color requires a Hasse arm of length $2(2t-1)$ and its midpoint. If the common arm has this length, it contributes one occurrence to both objects. The short varying arm contributes one additional occurrence only to $K_{t,M}$. This is the first Hasse histogram difference.

For graded rounds $s<t$, equality follows directly from Corollary~\ref{cor:object-sandwich}: the Hasse histograms agree at round $2s<2t$, and therefore the graded histograms agree at round $s$.

To show distinction at round $t$, use the auxiliary line-graph refinement in Lemma~\ref{lem:line-profile}. An edge on an arm of length at least two has line-graph degree $3$ at the center end, $1$ at the leaf end, and $2$ in between. For $t\geq2$, the midpoint edge $e_t$ of the short arm has degree $2$ and two stopped rays whose first non-$2$ marks are respectively $3$ and $1$, both at distance $t-1$. A radius-$(t-1)$ marked unfolding therefore identifies exactly this midpoint of an arm of length $2t-1$. Arms of lengths $2t$ and $2t+1$ contribute none. The common arm contributes one to both objects if its length is $2t-1$, and otherwise none; if $M=1$, its sole edge has degree $2$ with a $3$ mark in both directions and cannot match this profile. Thus the auxiliary histograms differ, and Lemma~\ref{lem:line-profile} makes the GWL histograms differ at round $t$.

When $t=1$, the edge loads on the varying arms are $\{2,3,2,1\}$ for lengths $(1,3)$ and $\{3,1,3,1\}$ for lengths $(2,2)$, so the first difference again occurs at round $t$. Since each new hash retains the old color, both histogram differences persist at all later rounds.
\end{proof}

\section{A One-Round Graph Criterion}

The following criterion characterizes equality of the two cell partitions after one round on a connected graph.

\begin{proposition}[One-round graph criterion]
\label{prop:graph-criterion}
Let $G$ be a connected simple graph, viewed as a one-dimensional complex with rank-only initialization. Its one-round Hasse--WL and GWL cell partitions agree if and only if $\deg u+\deg v$ is constant over $uv\in E(G)$. For a connected graph, this is equivalent to regularity or bipartite biregularity. Under this condition, the common one-round partition is stable for both refinement processes.
\end{proposition}

\begin{proof}
The first Hasse round colors each vertex by its degree and gives every edge cell the same color. GWL gives the same vertex partition. For an edge $uv$, its additional lower-mediated multiplicity is
\[
                 (\deg u-1)+(\deg v-1)=\deg u+\deg v-2.
\]
Thus the edge partitions agree exactly when the endpoint-degree sum is constant. If the constant is $D$, degrees along every path alternate as $a,D-a,a,D-a,\ldots$. An odd cycle forces $a=D-a$, and connectedness then forces regularity. If there is no odd cycle, the graph is bipartite with constant degree on each side, hence biregular. The converse is immediate. Under the condition, the vertex-degree classes and the single edge class have constant direct and mediated profiles, so the partition is equitable and stable.
\end{proof}

The path $P_4$, whose edge endpoint-degree sums are $3,4,3$, is a minimal example with different one-round partitions.

\clearpage
\begingroup
\small
\bibliographystyle{plain}
\bibliography{paper}
\endgroup

\end{document}

%% file: experiment_section.tex
\subsection{Depth and resource comparison}
\label{sec:depth-resource-experiment}
We compare a nonlinear four-channel layer with its two-round Hasse
factorization, using identical inputs and weights. Let
$\phi(h)=\tanh(hW_p)$ and $\psi(h)=\tanh(hW_w)$. With incidence signs
$\varepsilon_{\rho\sigma}$, the lower-mediated channel is
\[
L_\sigma=\sum_{\rho\lessdot\sigma}\varepsilon_{\rho\sigma}\psi(h_\rho)
 \odot\!\sum_{\substack{\rho\lessdot\sigma'\\\sigma'\ne\sigma}}
 \varepsilon_{\rho\sigma'}\phi(h_{\sigma'}),
\]
and $U_\sigma$ is the corresponding sum through cofaces. The direct
channels $B_\sigma,C_\sigma$ are signed sums of $\phi$. The update is
\[
h'_\sigma=\tanh\!\left(h_\sigma W_s+b_{\rank(\sigma)}
+\frac{0.2}{d}(B_\sigma+C_\sigma)
+\frac{0.3}{d^2}\bigl[1+0.25\tanh(h_\sigma W_g)\bigr]
 \odot(L_\sigma+U_\sigma)\right),
\]
where the Hasse-degree bound $d$ is supplied as a fixed model parameter.
The first Hasse round retains the input $h$ and two signed $q$-vector
star sums. The second removes one target contribution from each relevant
witness sum and applies the update. The intermediate state has $3q$
numeric coordinates and a separately exposed rank label.

The complexes are square grids with $k\times k$ vertices and a fixed
diagonal in each square, for $k=6,12,24$, and cones on cycles with
$16,64,128$ rim vertices. The grid cell counts are $171,771,3267$.
We use widths $q=4,16,64$, depths $T=1,2,4$, and ten independently
seeded inputs and layer sequences, giving 540 paired cases.
The largest coordinate error between the direct $T$-layer output and
the Hasse $2T$-round output is $1.67\times10^{-16}$.

For the shorter budget, let $\widetilde h_T$ be the stored cell signal
after $T$ Hasse rounds. It has completed $\lfloor T/2\rfloor$ layers;
on an odd round that signal is retained while the next sums are pending.
Table~\ref{tab:depth-resource-experiment} reports
$\varepsilon_T=\|\widetilde h_T-h^{\mathrm{direct}}_T\|_2/
\|h^{\mathrm{direct}}_T\|_2$, with norms over all cells and coordinates,
together with execution time and peak dynamic allocation.

\begin{table}[ht]
\centering\small
\begin{tabular}{lrrrrr}
\toprule
Family & Cells & Peer records & Time (ms) & Peak (MiB) & $\varepsilon_T$\\
\midrule
Grid & 171 & 1288 & 7.72 / 5.84 & 1.28 / 1.19 & 1.827 \\
Grid & 3267 & 28288 & 143.56 / 103.43 & 26.28 / 22.21 & 1.818 \\
Wheel & 65 & 528 & 2.55 / 1.81 & 0.70 / 0.47 & 1.821 \\
Wheel & 513 & 18560 & 65.50 / 14.80 & 18.19 / 3.32 & 1.809 \\
\bottomrule
\end{tabular}
\caption{Selected $q=64,T=4$ cases, with medians across ten seeds.
Time and peak allocation are shown as direct $T$ layers / Hasse $2T$
rounds. The direct reference materializes both mediated channels;
the Hasse implementation factors their sums. Allocation is measured
during execution, excluding prebuilt topology, inputs and weights.
$\varepsilon_T$ is the error of the partially completed $T$-round schedule.}
\label{tab:depth-resource-experiment}
\end{table}

At a wheel hub, explicit peer materialization creates quadratically
many records. Cover-only routing stores the incidences and fixed-width
sums instead. This explains the growing allocation difference in
Table~\ref{tab:depth-resource-experiment} and Figure~\ref{fig:depth-resources}.
The grid computations also show that twice as many communication rounds
can have lower execution cost when their messages are simpler.

The path lower-bound check uses $\Delta_0^T$ on a path of $T$ graph
edges. A left-endpoint impulse gives $(-1)^T$ at the opposite endpoint;
the zero input gives zero. Both inputs agree on every target Hasse ball
of radius less than $2T$. Executions at $T=1,2,3,4,6,8$ recover the
endpoint response first at round $2T$. Separately, joint rank-initialized
WL runs on 18 three-arm-tree pairs, with $t=1,2,3,4,6,8$ and
$M=1,2,9$, distinguish first at graded round $t$ and Hasse round $2t$.